\documentclass[12pt]{article}
\usepackage{pifont}

\usepackage{color}

\usepackage{amssymb}
\usepackage{amsfonts}
\usepackage{graphicx}
\usepackage{amsmath}
\usepackage{mathrsfs}
\usepackage{psfrag}
\usepackage{latexsym}
\usepackage{setspace}
\usepackage{amsbsy}
\usepackage{epsfig}
\usepackage{enumerate}
\usepackage{bm}

\newtheorem{thm}{Theorem}[section]

\newtheorem{coro}{Corollary}[section]

\newtheorem{lem}{Lemma}[section]
\newtheorem{definition}{Definition}[section]
\newtheorem{remark}{Remark}[section]

\newtheorem{assumption}{Assumption}
\newtheorem{conjecture}{Conjecture}
\newcommand{\Proof}{\textbf{Proof}}

\title{Almost Sure Central Limit Theorem in Sub-linear Expectation Spaces}

\date{}

\author{Weihuan Huang, Panyu Wu\\
{\small Zhongtai Securities Institute for Financial Studies, Shandong University}\\ {\small Jinan 250100, China; Email: huangweihuan@126.com; wupanyu@sdu.edu.cn}}

\begin{document}

\maketitle

\begin{abstract}
Peng (2006) initiated a new kind of central limit theorem under sub-linear expectations. Song (2017) gave an  estimate of the rate of convergence of Peng's central limit theorem. Based on these results, we %study the convergence rate of the central limit theorem, and
  establish a new kind of almost sure central limit theorem under sub-linear expectations in this paper, which is a quasi sure convergence version of Peng's central limit theorem. Moreover, this result is a natural extension of the classical almost sure central limit theorem to the case where the probability is no longer additive.
  %, and can be considered as a simulation method to some contingent claims with distribution uncertainty.
  Meanwhile, we prove a new kind of strong law of large numbers for non-additive probabilities without the independent identically distributed assumption.\\
\par  $\textit{Keywords:}$ Central limit theorem, law of large numbers, non-additive probability, sub-linear expectation.
\end{abstract}

%% \linenumbers

%% main text

\section{Introduction}
In recent years many interesting extensions of classical probability limit theorems involving log average and log density have been obtained. The basic result and starting point of these investigations is the almost sure central limit theorem established by Brosamler \cite{Brosamler} and Schatte \cite{Schatte}.
%\textit{
%Let $\{X_i\}_{i=1}^{\infty}$ a sequence of independent identically distributed random variables on a probability space $(\Omega, \mathcal{F}, P)$ with $EX_{1}=0$, $EX_{1}^2=1$ and $EX_{1}^{2+\alpha}<\infty$ for some $\alpha>0$, then
%\begin{equation}\label{ASCLT}
%\lim_{n\rightarrow \infty} \frac{1}{\log n}\sum_{k=1}^{n}\frac{1}{k}I\left[\frac{S_{k}}{\sqrt{k}}\in A\right]=\frac{1}{\sqrt{2\pi}}\int_{A}e^{-\frac{t^2}{2}}dt, \ a.s.,
%\end{equation}
%for any Borel-set $A\subset\mathbb{R}$ with $\lambda(\partial A)=0$, where $S_{k}:=\sum_{i=1}^{k}X_{i}$, $I$ denotes the indicator function and $\lambda$ is the Lebesgue measure.
%}
Let $\delta(x)(\cdot)$ denote the Dirac measure for $x\in\mathbb{R}$
%, i.e. $\delta(x)(A)=1$, if $x\in A$ and $=0$, if otherwise, for all continuity sets $A$, and let
and $N(0,1)(\cdot)$ denote the standard normal distribution on $\mathbb{R}$. Then the following result is the Theorem 1.2 in \cite{Brosamler}:\\
\textit{
Let $\{X_n\}_{n=1}^{\infty}$ be a sequence of independent identically distributed random variables on a probability space $(\Omega, \mathcal{F}, P)$ with $EX_{1}=0$, $EX_{1}^2=1$ and $E|X_{1}|^{2+2\alpha}<\infty$ for some $\alpha>0$, then $P$-a.s.
%\begin{equation}\label{ASCLT}
%\lim_{n\rightarrow \infty} \frac{1}{\log n}\sum_{k=1}^{n}\frac{1}{k}I\left[\frac{S_{k}}{\sqrt{k}}\in A\right]=\frac{1}{\sqrt{2\pi}}\int_{A}e^{-\frac{t^2}{2}}dt, \ a.s.,
%\end{equation}
\begin{equation}\label{1.1}
\lim_{n\to\infty}\frac{1}{\log n}\sum_{k=1}^{n}\frac{1}{k}\delta\left(\frac{S_{k}}{\sqrt{k}}\right)= N(0,1),
\end{equation}
where  $S_k=\sum_{i=1}^{k}X_i$ and the convergence is weak convergence of measures on $\mathbb{R}$.
%for any Borel-set $A\subset\mathbb{R}$ with $\lambda(\partial A)=0$, where $S_{k}:=\sum_{i=1}^{k}X_{i}$, $I$ denotes the indicator function and $\lambda$ is the Lebesgue measure.
}

Because of $$\int f d\left(\frac{1}{\log n}\sum_{k=1}^{n}\frac{1}{k}\delta\left(\frac{S_k}{\sqrt{k}}\right)\right) = \frac{1}{\log n}\sum_{k=1}^{n}\frac{1}{k}f\left(\frac{S_k}{\sqrt{k}}\right),$$
by Portmanteau Theorem (see Theorem 2.1 in \cite{C.P.M.}), (\ref{1.1}) is equivalent to that for any bounded and continuous function $f$, we have
\begin{equation}\label{1.2}
\lim_{n\rightarrow\infty}\frac{1}{\log n}\sum_{k=1}^{n}\frac{1}{k}f\left(\frac{S_{k}}{\sqrt{k}}\right)= \int fdN(0,1).
\end{equation}
Combining above equivalent relation, it follows from Theorem 11.3.3 in \cite{Dudley} that (\ref{1.1}) is equivalent to  (\ref{1.2}) holds for all bounded and Lipschitz continuous functions.

Later, Lacey and Philipp \cite{ASCLT1990} proved the almost sure central limit theorem only assume that $\{X_{n}\}_{n=1}^\infty$ are independent identically distributed with $EX_{1}=0$, $EX_{1}^2=1$. The use of the logarithmic mean looks perhaps peculiar, but it is essentially the only summation method that works (see Remark (d) in \cite{ASCLT1990}). For example, arithmetic means are not suitable for making the sequence $\left\{\delta\left(\frac{S_{n}}{\sqrt{n}}\right)\right\}_{n=1}^{\infty}$ weakly convergent to the standard normal distribution almost surely (see Theorem 1 in \cite{Schatte}).
Berkes and Dehling \cite{BH} revealed that if $\{X_n\}_{n=1}^{\infty}$ are independent identically distributed, then the almost sure central limit theorem is equivalent to the corresponding ordinary central limit theorem.
%$\frac{S_n}{\sqrt{n}}\rightarrow^{d} N(0,1)$, where $N(0,1)$ is the standard normal distribution.
However, it is different for some general situations, see Berkes \cite{B}, Berkes and Cs\'{a}ki \cite{BC}, Berkes and Dehling \cite{log}, \cite{BH}. %Cs\'{a}ki et al \cite{CFR}
These papers showed that every weak limit theorem (not only the central limit theorem) for independent random variables has an analogous almost sure version.

The key of these classical results are the theory of weak convergence and the additivity of the probabilities and expectations. However, such additivity assumption is not reasonable in many areas of applications, especially in economy and finance, because many uncertain phenomena can not be well modeled using additive probabilities or additive expectations, see for example, Artzner et al. \cite{jinrong3}, Avellaneda et al. \cite{jinrong1},  El Karoui et al. \cite{Peng1997}, Gilboa \cite{Gilboa} and Lyons \cite{jinrong2}. There are several different theories to deal with the uncertainty, such as imprecise probability (see \cite{Walley1991}), capacity theory (see \cite{Gilboa2}), possibility theory (see \cite{Dubois}) and fuzzy measure theory (see \cite{WangKlir}). In this paper, we use the framework of sub-linear expectation introduced by Peng \cite{Peng1}-\cite{Peng's book}.
Motivated by the risk measures, superhedge pricing and modelling uncertainty in finance, Peng initiated a lot of studies of independent identically distributed random variables, $G$-normal distribution and $G$-Brownian motion  under sub-linear expectations (see \cite{Peng1}-\cite{Peng's book}). These notions may be far-reaching in the sense that it is not based on a classical probability space given a priori. %, for example???, see???, we present their definitions in {\bf section 2}.

%Intuitively speaking, for a given set $\mathcal{P}$ of probability measures on $(\Omega,\mathcal{F})$, we define a pair of non-additive probabilities $(\mathbb{V}, \nu)$
% by
%$$\mathbb{V}(A)=\sup_{P\in \mathcal{P}}P(A),\quad  \nu(A)=\inf_{P\in\mathcal{P}}P(A),\quad \hbox{for all } A\in\mathcal{F},$$
%which are called upper probability and lower probability respectively.
%Let $(\mathbb{E},\mathcal{E})$ denote the upper-lower expectations generated by $\mathcal{P}$, that is,
%$$\mathbb{E}[X]=\sup_{P\in \mathcal{P}}E_{P}(X),\quad  \mathcal{E}[X]=\inf_{P\in\mathcal{P}}E_{P}(X).$$
%The set $\mathcal{P}$ represents the ambiguity of probability measures and $\mathbb{E}$ is a sub-linear expectation which is equivalent to the notion of coherent risk measures in \cite{jinrong3}.
%%$\nu(A)=1$ means that for all $P\in\mathcal{P}$, $P(A)=1$, in this case, we call $A$ holds quasi-surely. We would give an introduction of non-additive probability (sub-linear expectation) in {\bf section 2}.

One of the most important results given by Peng is the following central limit theorem under the sub-linear expectation which is a simplified version of Theorem 3.5 (the moment condition due to Remark 3.8) in Chapter II of \cite{Peng's book}:\\
\textit{
Let $\{X_{n}\}_{n=1}^{\infty}$ be a sequence of independent identically distributed random variables under sub-linear expectation $\mathbb{E}$, with $\mathbb{E}[X_1]=\mathbb{E}[-X_1]=0$
and $\mathbb{E}[|X_1^3|]<\infty$. Set $S_n=\sum_{i=1}^{n}X_i$.
%, $\mathbb{E}[X_1^2]=\overline{\sigma}^2$, $\mathcal{E}[X_1^2]=\underline{\sigma}^2>0$, and there exists a constant $\alpha>0$ such that $\mathbb{E}[|X_{1}|^{2+\alpha}]< \infty$.
Then for all bounded and Lipschitz continuous function $f:\mathbb{R}\to \mathbb R$,
\begin{equation*}
\lim_{n\rightarrow\infty}\mathbb{E}\left[f\left(\frac{S_n}{\sqrt{n}}\right)\right]= \widetilde{\mathbb E}\left[f(\xi)\right],
\end{equation*}
where $\xi$ is a $G$-normal distributed random variable under sub-linear expectation $\tilde{\mathbb{E}}$
%with variance $[\underline{\sigma}^2, \overline{\sigma}^2]$.
with $G(a)=\frac{1}{2}\mathbb{E}[aX_1^2]$, $a\in\mathbb{R}$.
}\\
Peng proved it by using the unique viscosity solution of a non-linear parabolic partial differential equation as the distribution of $G$-normal random variables. Recently, Song \cite{Song2018} gave an  estimate of the rate of convergence of Peng's central limit theorem for an independent but not necessarily identically distributed sequence of random variables under sub-linear expectation.
%And we will adopt Song's framework in our paper.

\iffalse
\textit{ Let $X_i$ be independent random variables under a sublinear expectation $\mathbb{E}$. We suppose that, for each $1\leq i\leq n$,  $\mathbb{E}[X_i]=\mathcal{E}[X_i]=0$, $\mathbb{E}[X_i^2]=\overline{\sigma}_i^2$, $\mathcal{E}[X_i^2]=\underline{\sigma}_i^2$, and denote $\sigma_i:= \frac{\overline{\sigma}_i + \underline{\sigma}_i}{2}$ and $B_n^2:= \sum_{i=1}^n \sigma_i^2$, then there exists $\alpha\in(0,1)$ depending on $\beta$, and $C_{\alpha,\beta}>0$ such that
$$
\underset{|f|_{Lip}\leq 1}\sup \bigg| \mathbb{E}\left[ f\left( \frac{S_n}{\sqrt{n} B_n} \right) \right] - \tilde{\mathbb{E}}[f(\xi)] \bigg| \leq C_{\alpha,\beta} \underset{1\leq i\leq n}\sup \bigg\{ \frac{\mathbb{E}[|X_i|^{2+\alpha}]}{\sigma_i^{2+\alpha}} \left( \frac{\sigma_i}{\sigma} \right)^\alpha \bigg\},
$$
where $\xi$ is $G$-normal distribution with $\beta:= \frac{\overline{\sigma}}{\underline{\sigma}}\geq 1$ and $\underline{\sigma} + \overline{\sigma}=2$.
}\fi

Hence, it is natural to think whether it is possible to have an analogous almost sure version of Peng's central limit theorem?
%\begin{problem}
%Is it possible to have an analogous quasi-sure version of Peng's CLT?
%\end{problem}
In this paper, we will find function $f$ such that (\ref{Q1}) holds%, %when $\{X_i\}_{i=1}^{\infty}$ are independent identically distributed under $\mathbb{E}$ with $\mathbb{E}[X_1]=\mathcal{E}[X_1]=0$ and $\mathbb{E}[X_1^2]=\overline{\sigma}^2$, $\mathcal{E}[X_1^2]=\underline{\sigma}^2$,
\begin{equation}\label{Q1}
\nu\left(\lim_{n\rightarrow \infty}\frac{1}{\log n}\sum_{k\leq n}\frac{1}{k}f\left(\frac{S_{k}}{\sqrt{k}}\right) = \widetilde{\mathbb{E}}\left[f(\xi)\right] \right)=1,
\end{equation}
where $\nu$ is the lower probability generated by the sub-linear expectation $\mathbb{E}$.
%We give (\ref{Q1}) a name the quasi sure central limit theorem.

Motivated by Lacey and Philipp \cite{ASCLT1990}, we will prove (\ref{Q1}) through the following strong law of large numbers,
\begin{equation}\label{Q3}
\nu\left( \limsup_{N\rightarrow \infty} \frac{1}{N}\sum_{l=1}^{N} Z_{l}\leq 0 \right)=1,
\end{equation}
where $Z_{l}:=\sum_{4^{l-1}\leq k < 4^{l}} \frac{\xi_{k}}{k}$, $\xi_{k}:=f(\frac{S_{k}}{\sqrt{k}})-\mathbb{E}f(\frac{S_{k}}{\sqrt{k}})$.
%We put the details in section 3.
Chen \cite{Chen4} established a strong law of large numbers
\begin{equation}\label{chen'sSLLN}
\nu\left( \underline{\mu}\leq \liminf_{n\rightarrow \infty}\frac{\sum_{l=1}^{n}Z_{l}}{n} \leq \limsup_{n\rightarrow \infty}\frac{\sum_{l=1}^{n}Z_{l}}{n}\leq \overline{\mu}\right)=1,
\end{equation}
when $\{Z_l\}_{l=1}^{\infty}$ are independent identically distributed under $\mathbb{E}$ with $\mathbb{E}[Z_1]=\overline{\mu}$ and $-\mathbb{E}[-Z_1]=\underline{\mu}$. Chen et al. \cite{CHW} established an extension strong law of large numbers under exponential independence.
Unfortunately, the $\{Z_l\}_{l=1}^{\infty}$ in (\ref{Q3}) are  neither independent  nor exponential independent under $\mathbb{E}$, thus we give a strong law of large numbers without independence but with some moment conditions in section 3.
 %\\
%\textbf{Theorem}   \textit{
   %Let $T_{n}:=\sum_{l=1}^{n}Z_{l}$. If  for any $ n \geq 1$ and arbitrary $ n_{1}\le n_{2},$ we have
   %\begin{equation*}
   %\begin{cases}
   %\mathbb{E}[T_{n}]\leq 0, \\
   %\mathbb{E}[T_{n}^{2}]\leq M_{0}\cdot n + [\mathbb{E}T_{n}]^{2}, \ (or \ \mathbb{E}[T_{n}^{2}]\leq M_{0}\cdot n ), \\
   %\sum_{l=n_{1}}^{n_{2}}\mathbb{E}(Z_{l}^2)\leq M_{1}\cdot (n_{2}-n_{1}+1),
   %\end{cases}
   %\end{equation*}
   %where $M_{0}$ and $M_{1}$ are constants, then
   %\begin{equation*}
   %\nu\left( \limsup_{n\rightarrow \infty}\frac{T_{n}}{n}\leq 0\right)=1.
   %\end{equation*}}
%We put this Theorem and some corollaries in appendix B.
This strong law of large numbers drives us to verify whether the $\{Z_l\}_{l=1}^{\infty}$ in (\ref{Q3}) satisfy the moment conditions. So, we restrict the test functions $f$ in (\ref{Q1}) to be elements in linear space $\mathcal{H}=\{f\in C_{b,Lip}(\mathbb{R}):\widetilde{\mathbb{E}}[f(\xi)]=-\widetilde{\mathbb{E}}[-f(\xi)]]\}$. %, expectation invariant functions of $G$-normal distribution.
%And we need to study convergence rate of Peng's central limit theorem in appendix A to get $\mathbb{E}[T_{n}^{2}]\leq M_{0}\cdot n$.
%In {\bf section 3}, we show that all the odd functions and all their linear translations are in $\tilde{\mathcal{H}}$. %furthermore for any continuous $G$-convex function $f$,
%%\begin{equation}\label{Jensen}
%%\nu\left(\lim_{n\rightarrow \infty}f\left(\frac{1}{\log n}\sum_{k\leq n}\frac{1}{k}\cdot \frac{S_k}{\sqrt{k}} \right) \leq \widetilde{\mathbb{E}}\left[f(\xi)\right] \right)=1.
%%\end{equation}
%We also show that if the non-additive probability $\nu$ in (\ref{Q1}) is a probability measure, then (\ref{Q1}) becomes \eqref{ASCLT}, hence the quasi sure central limit theorem (almost sure central limit theorem for non-additive probabilities) is a natural extension of the classical almost sure central limit theorem.
%And we have \eqref{Q1} holds for all bounded Lipschitz functions when $\xi$ is a classical normal distribution random variable in a sublinear expectation space.
Hu \cite{Hu} considered different moment conditions to prove the law of large numbers in sub-linear expectation spaces.
This paper is organized as follows. In section 2, we introduce some basic concepts and lemmas in the sub-linear expectation theory which will be used in the sequel. In section 3, we give a new kind of strong law of large numbers for non-additive probabilities under some moment conditions. In section 4, we investigate the almost sure central limit theorem in sub-linear expectation spaces.

\section{Preliminaries}

We adopt the framework and notations in Peng \cite{Peng5} and \cite{Peng's book}.
%and Song \cite{Song2018}.
Let $(\Omega,\mathcal F)$ be a given measurable space and $\mathscr{H}$ be a linear space of real measurable functions defined on $(\Omega,\mathcal F)$ such that if $X_1,\ldots, X_n \in \mathscr{H}$  then $f(X_1,\ldots,X_n)\in \mathscr{H}$, for each $f\in C_{b,Lip}(\mathbb R^n)$, where $C_{b,Lip}(\mathbb R^n)$ denotes the space of all bounded and Lipschitz continuous functions on $\mathbb R^n$.
%for each $f\in C_b(\mathbb R^n)\bigcup  C_{l,Lip}(\mathbb R^n)$,  where $C_b(\mathbb R^n)$ denotes the space of all bounded continuous functions and $C_{l,Lip}(\mathbb R^n)$ denotes the linear space of local Lipschitz functions $f$ satisfying
%\begin{eqnarray*} & |f(\bm x) - f(\bm y)| \le  C(1 + |\bm x|^m + |\bm y|^m)|\bm x- \bm y|, \;\; \forall \bm x, \bm y \in \mathbb R^n,&\\
%& \text {for some }  C > 0, m \in \mathbb  N \text{ depending on } f. &
%\end{eqnarray*}
%Further, we let $C_{b,Lip}(\mathbb R^n)$ denote the space of all bounded and Lipschitz functions on $\mathbb R^n$.
\begin{definition}
A {\bf sub-linear expectation} $\mathbb{E}$ on $\mathscr{H}$  is a functional $\mathbb{E}: \mathscr{H}\to \mathbb R$ satisfying the following properties: for all $X, Y \in \mathscr H$, we have
\begin{description}
\item[\rm (a)] {\bf Monotonicity}: $\mathbb{E} [X]\ge \mathbb{E} [Y]$, if $X \ge  Y$;
\item[\rm (b)] {\bf Constant preserving}: $\mathbb{E} [c] = c$, for $c \in \mathbb{R}$;
\item[\rm (c)] {\bf Sub-additivity}: $\mathbb{E}[X+Y]\le \mathbb{E} [X] +\mathbb{E} [Y]$;
\item[\rm (d)] {\bf Positive homogeneity}: $\mathbb{E} [\lambda X] = \lambda \mathbb{E}[X]$,  for $ \lambda\ge 0$.
%\item[\rm (e)] $\mathbb{E}[f_n(X)]\downarrow 0$, for $X\in\mathscr{H}$ and $f_n\in C_{b,Lip}(\mathbb R^n)$, $f_n\downarrow 0$.
 \end{description}
 The triple $(\Omega, \mathscr{H}, \mathbb{E})$ is called a sub-linear expectation space.
 \end{definition}
%For $X\in\mathscr{H}$, set
%$$
%\mathcal{N}^X[f]:= \mathbb{E}[f(X)], \ f\in C_{b,Lip}(\mathbb R^n),
%$$
%which is also a sublinear expectation on $C_{b,Lip}(\mathbb R^n)$. We say $X$ is distributed as $\mathcal{N}^X$. A functional $\mathcal{N}$ is a sublinear expectation on $C_{b,Lip}(\mathbb R^n)$ if and only if it can be represented as the supremum expectation of a weakly compact subset $\Theta$ of probability measures on $(\mathbb{R}, \mathcal{B}(\mathbb{R}))$, (see \cite{Denis}),
%$$
%\mathcal{N}[f] = \underset{\mu\in\Theta}\sup \mu[f], \;\; \hbox{for all } f\in C_{b,Lip}(\mathbb R^n).
%$$

By Lemma 2.4 in  Peng \cite{Peng5},  a functional $\mathbb{E}$ defined on $(\Omega,\mathscr{H})$ is a sub-linear expectation if and only if there exists a family of linear expectations (finite additive) $\{ E_\theta : \theta\in\Theta\}$ such that  $\mathbb{E}[X]=\sup_{\theta\in\Theta}E_{\theta}[X]$,  for all $X\in\mathscr{H}.$
%Let us denote the conjugate expectation $\mathcal{E}$ of $\mathbb{E}$ by
%$ \mathcal{E}[X]=-\mathbb{E}[-X]$,  for all $X\in \mathscr{H}. $
In the sequel, we will consider each $E_\theta$ is generated by a $\sigma$-additive probability measure,
that is
$$\mathbb{E}[X]=\underset{{P\in \mathcal{P}}}\sup E_{P}[X],\ \ \hbox{for all } X\in \mathscr{H},$$
%\ \  \mathcal{E}[X]=\underset{{P\in \mathcal{P}}}\inf E_{P}[X],$$
where $\mathcal{P}$ is a set of $\sigma$-additive probability measures and $E_P$ is the expectation with respect to $P$.
Then, we define the upper probability $\mathbb{V}$ and lower probability $\nu$ by
$$\mathbb{V}(A)=\sup_{P\in\mathcal{P}}P(A), \; \nu(A)=\inf_{P\in\mathcal{P}}P(A),\;\; \hbox{for all } A\in\mathcal{F}.$$
It is easy to check that $\mathbb{V}$ is continuous from below and $\nu$ is continuous from above, in other words, $\mathbb{V}(A_{n})\uparrow \mathbb{V}(A)$ if $A_{n}\uparrow A$, and $\nu(A_{n})\downarrow \nu(A)$ if $A_{n}\downarrow A$, where $A_{n}, A\in {\cal F}, n\ge 1$.
\begin{definition}
A set $A$ is a polar set if $\mathbb{V}(A)=0$ and a property holds quasi-surely if it holds outside a polar set.
\end{definition}
%We call almost sure central limit theorem for non-additive probabilities quasi sure central limit theorem for the reason that the limit in (\ref{Q1}) %and (\ref{Q2})
%holds quasi-surely.
\begin{definition}
{\rm (i)} Let $\bm X_1$ and $\bm X_2$ be two $n$-dimensional random vectors defined
respectively in sub-linear expectation spaces $(\Omega_1, \mathscr{H}_1, \mathbb{E}_1)$
  and $(\Omega_2, \mathscr{H}_2, \mathbb{E}_2)$. They are called {\bf identically distributed}, denoted by $\bm X_1\overset{d}= \bm X_2$  if
$$ \mathbb{E}_1[f(\bm X_1)]=\mathbb{E}_2[f(\bm X_2)], \ \hbox{for all } f \in C_{b,Lip}(\mathbb R^n).$$
%A sequence $\{X_n\}_{n=1}^{\infty}$ of random variables is said to be identically distributed if $X_n\overset{d}= X_1$ for each $n\ge 1$.

{\rm (ii)}  In a sub-linear expectation space  $(\Omega, \mathscr{H}, \mathbb{E})$, a random vector $\bm Y \in\mathscr{H}^n$ is said to be {\bf independent} to a random vector $\bm X
\in \mathscr{H}^m$ under $\mathbb{E}$  if for each test function $f \in C_{b,Lip}(\mathbb R^{m+n})$
we have
$$ \mathbb{E} [f(\bm X, \bm Y )] = \mathbb{E} \big[\mathbb{E}[f(\bm x, \bm Y )]\big|_{\bm x=\bm X}\big],$$
whenever $\overline{f}(\bm x):=\mathbb{E}\left[|f(\bm x, \bm Y )|\right]<\infty$ for all $\bm x$ and
 $\mathbb{E}\left[|\overline{f}(\bm X)|\right]<\infty$.

{\rm (iii)}  A sequence of random variables $\{X_n\}_{n=1}^{\infty}$
 is said to be {\bf independent}, if $X_{n+1}$ is independent to $(X_1,\ldots, X_n)$ for each $n\ge 1$.
\end{definition}

Next, we recall four lemmas which will be used in the following two sections.

\iffalse
\begin{lem}[Central Limit Theorem]\label{CLT1} (Theorem 3.5 and Remark 3.8 in Chapter II of Peng \cite{Peng's book}).
Suppose $\{X_{i}\}_{i=1}^{\infty}$ is a sequence of independent identically distributed random variables with $\mathbb{E}[X_1]=\mathcal{E}[X_1]=0$, $\mathbb{E}[X_1^2]=\overline{\sigma}^2$, $\mathcal{E}[X_1^2]=\underline{\sigma}^2$, and there exists a constant $\alpha>0$ such that $\mathbb{E}[|X_{1}|^{2+\alpha}]< \infty$. Then for all $f \in C_{b.lip}(\mathbb{R})$, we have
\begin{equation}\label{Peng's CLT}
\lim_{n\rightarrow\infty}\mathbb{E}\left[f\left(\frac{S_n}{\sqrt{n}}\right)\right]= \widetilde{\mathbb E}\left[f(\xi)\right].
 \end{equation}
\end{lem}
\fi
\iffalse
\begin{lem}[Donsker's Invariance Principle]\label{Donsker} (Theorem 3.1 in Zhang \cite{Zhang2015}).
Suppose $\{X_{i}\}_{i=1}^{\infty}$ is a sequence of independent identically distributed random variables with $\mathbb{E}[X_1]=\mathcal{E}[X_1]=0$, $\mathbb{E}[X_1^2]=\overline{\sigma}^2$ and $\mathcal{E}[X_1^2]=\underline{\sigma}^2$, and there exists a constant $\alpha>0$ such that $\mathbb{E}[|X_{1}|^{2+\alpha}]< \infty$. Then for all $f\in C_{b.lip}(C[0,1])$, we have
\begin{equation}\label{Zhang's CLT}
\lim_{n\rightarrow\infty}\mathbb{E}\left[f\left( W_{n} \right)\right] = \widetilde{\mathbb E}\left[f(W)\right].
 \end{equation}
\end{lem}
\fi

 \begin{lem}[Rosenthal's inequality]\label{Rosenthal} (Lemma 3.1 in Zhang \cite{Zhang2015}).
   Let $\{X_1,\ldots, X_n\}$ be a sequence  of  independent random variables
 in $(\Omega, \mathscr{H}, \mathbb{E})$. If $\mathbb{E}[X_k]=\mathbb{E}[-X_k]= 0$, $k=1,\cdots, n$, then for $p\ge 2$, we have
$$
\mathbb{E}\left[\max_{k\le n}|S_k|^p \right]\le \frac{C_p}{2}\left\{ \sum_{k=1}^n \mathbb{E} [|X_k|^p]+\left(\sum_{k=1}^n \mathbb{E} [|X_k|^2]\right)^{p/2}\right\},
$$
where $C_{p}$ is a positive constant depending only on $p$.
\end{lem}

\begin{lem}[Chebyshev's inequality] \label{Chebyshev} (Proposition 2.1(2) in Chen et al. \cite{Chen3}).
Let $g(x)>0$ be any given even function on $\mathbb{R}$ and nondecreasing on $(0,\infty)$. Then for any $x>0$,
$$\mathbb{V}(|X|\geq x)\leq \frac{\mathbb{E}[g(X)]}{g(x)}.$$
\end{lem}

\begin{lem}[Borel-Cantelli lemma] \label{Borel} (Lemma 2.2 in Chen et al. \cite{Chen3}).
If $\sum_{n=1}^{\infty}\mathbb{V}(A_{n})<\infty$, then
\begin{equation*}
\mathbb{V}\left(\bigcap_{k=1}^{\infty}\bigcup_{n=k}^{\infty}A_{n}\right)=0.
\end{equation*}
\end{lem}

\begin{lem}[H\"{o}lder's inequality] \label{Holder}(Proposition 2.1(1) in Chen et al. \cite{Chen3}).
For $p, q>1$ with $\frac{1}{p}+\frac{1}{q}=1$, we have
$$\mathbb{E}|XY|\leq (\mathbb{E}|X|^p)^{\frac{1}{p}}(\mathbb{E}|Y|^q)^{\frac{1}{q}}.$$
\end{lem}

%=================================================================================================================

\section{Strong law of large numbers for non-additive probabilities}

%Before give a kind of strong laws of large numbers (Theorem \ref{SLLN}), we first give a special case (Theorem \ref{thm2}), which has been used in the proof of quasi sure central limit theorem.

In order to distinguish the basic sequence $\{X_{n}\}_{n=1}^{\infty}$ in the almost sure central limit theorem, we let $\{Z_l\}_{l=1}^{\infty}$ denote a sequence of random variables in the sub-linear expectation $(\Omega,\mathscr{H},\mathbb{E})$. Set $T_{n}:=\sum_{l=1}^{n}Z_{l}$, for $n\ge 1$. We have the following law of large numbers under some general moment conditions.
% which is a special case of Theorem \ref{th2}.

\begin{thm}\label{thm2}
%\noindent\textbf{Theorem B.1.}
%\textit{
	For all $n \geq 1$ and  $n_{1}\leq n_{2},$ if
	\begin{equation*}
		\begin{cases}
			\mathbb{E}[T_{n}]\leq 0, \\
			\mathbb{E}[T_{n}^{2}]\leq M_{0}\cdot n + [\mathbb{E}T_{n}]^{2}, \ (or \ \mathbb{E}[T_{n}^{2}]\leq M_{0}\cdot n ), \\
			\sum_{l=n_{1}}^{n_{2}}\mathbb{E}(Z_{l}^2)\leq M_{1}\cdot (n_{2}-n_{1}+1),
		\end{cases}
	\end{equation*}
	where $M_{0}$ and $M_{1}$ are positive constants, then
	\begin{equation*}
		\nu\left( \limsup_{n\rightarrow \infty}\frac{T_{n}}{n}\leq 0\right)=1.
\end{equation*}%}
\end{thm}

\iffalse
\begin{remark}
	If $Z_{l}$ is i.i.d. and $\mathbb{E}[Z_{1}]=0$, $\mathbb{E}[Z_{1}^2]<\infty$, then $\mathbb{E}[T_{n}]\leq 0$, and $\mathbb{E}[T_{n}-\mathbb{E}[T_{n}]]^2=n\mathbb{E}[Z_{1}-\mathbb{E}Z_{1}]^2\leq Mn$, $\sum_{l=n_{1}}^{n_{2}}\mathbb{E}(Z_{l}^2)\leq M_{1}\cdot (n_{2}-n_{1}+1)$. Hence, the moment restrictions is weaker than the i.i.d. condition.
\end{remark}
\fi

\Proof.
Since $\mathbb{E}[T_{n}]\leq 0$ and $\mathbb{E}[T_{n}^{2}]\leq M_{0}\cdot n + [\mathbb{E}T_{n}]^{2}$, by the sub-additivity and positive homogeneity of $\mathbb{E}$, we have
\begin{align*}
	\mathbb{E}\left[(T_{n}-\mathbb{E}[T_{n}])^2\right] &= \mathbb{E}\left[T_{n}^2+(\mathbb{E}T_{n})^2-2T_{n}\mathbb{E}[T_{n}]\right] \\
	&\leq \mathbb{E}[T_{n}^2]+(\mathbb{E}T_{n})^2-2\mathbb{E}[T_{n}]\mathbb{E}[T_{n}] \\
	&= \mathbb{E}[T_{n}^2]-(\mathbb{E}T_{n})^2 \\
	&\leq M_{0}\cdot n.
\end{align*}
By Chebyshev's inequality (Lemma \ref{Chebyshev}), for any $ \epsilon > 0$, we have
\begin{equation*}
	\mathbb{V}\left( \left|\frac{T_{n}-\mathbb{E}T_{n}}{n}  \right| > \epsilon \right) \leq \frac{\mathbb{E}[(T_{n}-\mathbb{E}T_{n})^2]}{n^2\epsilon^2}
	\leq \frac{M_{0}}{n\epsilon^2}.
\end{equation*}
Hence, for the subsequence $\{n^2\}_{n=1}^{\infty}$,
\begin{equation*}
	\sum_{n=1}^{\infty}\mathbb{V}\left( \left|\frac{T_{n^2}-\mathbb{E}T_{n^2}}{n^2}  \right| > \epsilon \right) \leq \sum_{n=1}^{\infty}\frac{M_{0}}{n^2\epsilon^2}< \infty.
\end{equation*}
Due to Borel-Cantelli Lemma (Lemma \ref{Borel}), we have
\begin{equation*}
	\mathbb{V}\left( \bigcap_{k=1}^{\infty}\bigcup_{n=k}^{\infty}\left\{\left|\frac{T_{n^2}-\mathbb{E}T_{n^2}}{n^2}  \right| > \epsilon\right\}\right)=0.
\end{equation*}
Since $\mathbb{V}$ is continuous from below, we have
\begin{eqnarray*}
	\mathbb{V}\left( \limsup_{n \rightarrow \infty} \frac{T_{n^2}-\mathbb{E}T_{n^2}}{n^2}> 0 \right)
	&\le&\mathbb{V}\left( \bigcup_{\epsilon > 0}\bigcap_{k=1}^{\infty}\bigcup_{n=k}^{\infty}\left\{\left|\frac{T_{n^2}-\mathbb{E}T_{n^2}}{n^2}  \right| > \epsilon\right\}\right)\\
	&=& \lim_{\epsilon \downarrow 0}\mathbb{V}\left(\bigcap_{k=1}^{\infty}\bigcup_{n=k}^{\infty}\left\{\left|\frac{T_{n^2}-\mathbb{E}T_{n^2}}{n^2}  \right| > \epsilon\right\}\right)=
	0.
\end{eqnarray*}
%which is equivalent to
%\begin{equation*}
%\nu\left( \limsup_{n \rightarrow \infty} \frac{T_{n^2}-\mathbb{E}T_{n^2}}{n^2}\leq 0 \right)=1.
%\end{equation*}
Because of $\mathbb{E}[T_{n^2}]\leq 0$, we have
$$\left\{ \limsup_{n \rightarrow \infty} \frac{T_{n^2}}{n^2}> 0 \right\}\subseteq\left\{ \limsup_{n \rightarrow \infty} \frac{T_{n^2}-\mathbb{E}T_{n^2}}{n^2}> 0 \right\}.$$
Thus by the monotonicity of $\mathbb{V}$, we have
%\begin{equation*}
%\nu\left( \limsup_{n \rightarrow \infty} \frac{T_{n^2}}{n^2}\leq 0 \right)=1,
%\end{equation*}
%which is equivalent to
\begin{equation}\label{proofthm224}
	\mathbb{V}\left( \limsup_{n \rightarrow \infty} \frac{T_{n^2}}{n^2}> 0 \right)=0.
\end{equation}
Define $D_{n}:=\max_{n^2\leq k<(n+1)^2} \left| T_{k}- T_{n^2}\right|$, then
\begin{align*}
	\mathbb{E}[D_{n}^2] &=\mathbb{E}\left[\left(\max_{n^2\leq k<(n+1)^2} \left| T_{k}- T_{n^2}\right|\right)^2\right]
	=\mathbb{E}\left[\max_{n^2< k<(n+1)^2} \left| \sum_{l=n^2+1}^{k}Z_{l}\right|^2\right] \\
	&\leq \mathbb{E}\left[\max_{n^2< k<(n+1)^2} (k-n^2) \sum_{l=n^2+1}^{k}Z_{l}^2\right] \leq (2n+1)\mathbb{E}\left[ \sum_{l=n^2+1}^{(n+1)^2}Z_{l}^2 \right] \\
	&\leq (2n+1)\sum_{l=n^2+1}^{(n+1)^2}\mathbb{E}[ Z_{l}^2 ]
	\leq M_{1}(2n+1)(2n+1)
	\leq M_{2}\cdot n^2,
\end{align*}
where $M_{2}$ is a positive constant.
By Chebyshev's inequality (Lemma \ref{Chebyshev}), we have
\begin{equation*}
	\mathbb{V}(D_{n}>n^2\epsilon)\leq \frac{\mathbb{E}[D_{n}^2] }{n^4\epsilon^2}\leq \frac{M_{2}}{n^2\epsilon^2}.
\end{equation*}
Hence,
\begin{equation*}
	\sum_{n=1}^{\infty}\mathbb{V}(D_{n}>n^2\epsilon)\leq \sum_{n=1}^{\infty}\frac{M_{2}}{n^2\epsilon^2}<\infty.
\end{equation*}
By Borel-Cantelli Lemma (Lemma \ref{Borel}), we have
$$
\mathbb{V}\left(\bigcap_{k=1}^{\infty}\bigcup_{n=k}^{\infty} \left\{\frac{D_{n}}{n^2}>\epsilon\right\}\right)=0.
$$
It follows from the continuity from below and monotonicity of $\mathbb{V}$ that
\begin{eqnarray*}
	\mathbb{V}\left(\limsup_{n\rightarrow \infty} \frac{D_{n}}{n^2}>0\right)&=&\mathbb{V}\left(\bigcup_{\epsilon>0}\left\{\limsup_{n\rightarrow \infty} \frac{D_{n}}{n^2}>\epsilon\right\}\right)
	\\
	&=&\lim\limits_{\epsilon\downarrow0}\mathbb{V}\left(\limsup_{n\rightarrow \infty} \frac{D_{n}}{n^2}>\epsilon\right)\\
	&\leq &  \lim\limits_{\epsilon\downarrow0}\mathbb{V}\left(\bigcap_{k=1}^{\infty}\bigcup_{n=k}^{\infty} \left\{\frac{D_{n}}{n^2}>\epsilon\right\}\right)=0.
\end{eqnarray*}
That is
\begin{equation}\label{proofthm23}
	\mathbb{V}\left(\limsup_{n\rightarrow \infty} \frac{D_{n}}{n^2}>0\right)=0.
\end{equation}
For $n^2\leq k<(n+1)^2$, $$\frac{T_{k}}{k}\leq \frac{T_{n^2}}{n^2} + \frac{D_{n}}{n^2}.$$
We deduce from (\ref{proofthm224}) and (\ref{proofthm23}) that
\begin{equation*}
	\mathbb{V}\left( \limsup_{n \rightarrow \infty} \frac{T_{n}}{n}> 0 \right) \leq \mathbb{V}\left( \limsup_{n \rightarrow \infty} \frac{T_{n^2}}{n^2}> 0 \right) + \mathbb{V}\left(\limsup_{n\rightarrow \infty} \frac{D_{n}}{n^2}>0\right) =0,
\end{equation*}
thus
\begin{equation*}
	\nu\left( \limsup_{n\rightarrow \infty}\frac{T_{n}}{n}\leq 0\right)=1.
\end{equation*}
Therefore, the proof of Theorem \ref{thm2} is completed. $\Box$

\section{Almost sure central limit theorem in sub-linear expectation space}
In this section, we consider the almost sure central limit theorem for an independent but not necessarily identically distributed sequence of random variables.
%First of all, we introduce some notations.
For a random variable $X$ in a sub-linear expectation space $(\Omega,\mathscr{H}, \mathbb{E})$ with $\mathbb{E}[X]=\mathbb{E}[-X]=0$, set $\beta:
= \frac{\sqrt{\mathbb{E}[X^2]}}{\sqrt{-\mathbb{E}[-X^2]}}$ and $\sigma:=\frac{\sqrt{\mathbb{E}[X^2]} + \sqrt{-\mathbb{E}[-X^2]}}{2} $ to characterize the variances of a random variable $X$ under a sublinear expectation space. Throughout this section, the random variable sequence is supposed to satisfy the following two assumptions:
\begin{assumption}\label{assump0}
The sub-linear expectation spaces $(\Omega,\mathscr{H},\mathbb{E})$ and $(\widetilde{\Omega},\widetilde{\mathscr{H}},\widetilde{\mathbb{E}})$ considered in this section satisfy that for all  $X\in\mathscr{H}(\hbox{or }\widetilde{\mathscr{H}})$ and $f_n\in C_{b,Lip}(\mathbb R)$, $f_n\downarrow 0$ we have $\mathbb{E}(\hbox{or }\widetilde{\mathbb{E}})[f_n(X)]\downarrow 0$.
\end{assumption}

\begin{assumption}\label{assump}
Fixed the ratio $\beta\geq 1$ of variances as a constant, the sequence $\{X_i\}_{i=1}^{\infty}$ is a sequence of independent random variables in any fixed sub-linear expectation space $(\Omega, \mathscr{H}, \mathbb{E})$  with  $\mathbb{E}[X_i]=\mathbb{E}[-X_i]=0$, $\overline{\sigma}_i=\sqrt{\mathbb{E}[X_i^2]}$ and $\underline{\sigma}_i=\sqrt{-\mathbb{E}[-X_i^2]}$, for $i\geq 1$, and $0< \inf\limits_{i\geq 1}\underline{\sigma}_i^2 \leq \sup\limits_{i\geq 1}\overline{\sigma}_i^2 < \infty$.  %$G$-Brownian motion $W:=\{W(t)\}_{t\in[0,1]}$ we considered is the  canonical process  on $(\widetilde{\Omega}, \widetilde{\mathscr H}, \widetilde{\mathbb E})$  with $\xi \sim N(0,[\underline{\sigma}^2,\overline{\sigma}^2])$ denotes $G$-normal distribution with uncertainty variance $[\underline{\sigma}^2, \overline{\sigma}^2]$.
For $n\ge 1$, denote $S_0=0$, $S_n= \sum\limits_{i=1}^n X_i$ and $\sigma_n:= \frac{\underline{\sigma}_n +\overline{\sigma}_n}{2}$, $\beta_n:= \frac{\overline{\sigma}_n}{\underline{\sigma}_n} \equiv \beta$, $B_n:= \sqrt{\sum\limits_{i=1}^n \sigma_i^2}$, $W_n:= \frac{S_n}{B_n}$.
\end{assumption}

 A Random variable $\xi$ is {\bf $G$-normal distributed} (denoted by $\xi\sim N(0, [\underline{\sigma}^2, \overline{\sigma}^2])$) under a sub-linear expectation $\widetilde{\mathbb{E}}$, if and only if for any $f \in C_{b,Lip}(\mathbb R)$, the function $u(t,x)=\widetilde{\mathbb{E}}\left[f\left(x+\sqrt{t} \xi \right)\right]$ ($x\in \mathbb R, t\ge 0$) is the unique viscosity solution of the following $G$-heat equation
\begin{equation*}
		\begin{cases}
			 \partial_t u -G\left( \partial_{xx}^2 u\right) =0,\ (x,t)\in\mathbb{R}\times(0,\infty), \\
			u(0,x)=f(x),
		\end{cases}
	\end{equation*}
where $G(a)=\frac{1}{2}\widetilde{\mathbb{E}}[a\xi^2]$, $a\in\mathbb{R}$, is determined by the variances $\bar{\sigma}^2:=\widetilde{\mathbb{E}}[\xi^2]$ and $\underline{\sigma}^2:=-\widetilde{\mathbb{E}}[-\xi^2]$. If $\bar{\sigma}^2=\underline{\sigma}^2$, then $G$-normal distribution is just the classical normal distribution $N(0, \bar{\sigma}^2)$.

%In this paper, we use $\beta:
%= \frac{\mathbb{E}[X^2]^\frac12}{\mathcal{E}[X^2]^\frac12} := \frac{\overline{\sigma}}{\underline{\sigma}} $ and $\sigma:=\frac{\mathbb{E}[X^2]^\frac12 + \mathcal{E}[X^2]^\frac12}{2} := \frac{\overline{\sigma} + \underline{\sigma}}{2}$ to characterize the variances of a random variable $X$ under any sublinear expectation space. And we shall write $\mathcal{N}_\beta(0,\sigma^2)$ for the $G$-normal distribution $\mathcal{N}_G$, and write $\mathcal{N}_\beta$ for $\mathcal{N}_\beta(0,1)$. Clearly, $\mathcal{N}_1(0,\sigma^2)= N(0,\sigma^2)$, the classical normal distribution.
%Let
%$$
%   \mathcal{H}:=\{f\in C_{b,Lip}(\mathbb{R}):\ \widetilde{\mathbb{E}}[f(\xi)]=\widetilde{\mathcal{E}}[f(\xi)]\},
%   $$

Let
$$
   \mathcal{H}:=\{f\in C_{b,Lip}(\mathbb{R}):\ \widetilde{\mathbb{E}}[f(\xi)]=-\widetilde{\mathbb{E}}[-f(\xi)]\},
   $$
where $\xi$ is $G$-normal distribution under the sub-linear expectation $\widetilde{\mathbb{E}}$. It is easy to check that $\mathcal{H}$ is a linear space, that is, if $f_1, f_2 \in \mathcal{H}$, then $a f_1+b f_2 \in \mathcal{H}$, for all $ a,b \in \mathbb{R}$. %In this paper, we let the test functions $f$ in the ASCLT are the elements in $\mathcal{H}$ in order to calculus the covariance.
   In particular, if $f \in \mathcal{H}$, then $af+b \in \mathcal{H}$, for all $ a,b\in \mathbb{R}$. %Notice that $\mathcal{H}$ contains all odd

\iffalse
   \begin{thm}[Almost Sure Central Limit Theorem]%\label{thm1}
   Suppose that $\{X_{i}\}_{i=1}^{\infty}$ is a sequence of independent identically distributed random variables under sub-linear expectation $\mathbb{E}$ with $\mathbb{E}[X_{1}]=\mathcal{E}[X_{1}]=0, \ \mathbb{E}[X_{1}^2]=\overline{\sigma}^2, \ \mathcal{E}[X_{1}^2]=\underline{\sigma}^2$, and there exists a constant $\alpha>0$ such that $\mathbb{E}[|X_{1}|^{2+\alpha}]< \infty$.
   Then for all $f \in \mathcal{H}$,
   %$$\forall f \in \mathcal{H}=\{f\in C_{b,Lip}(C[0,1]):\widetilde{\mathbb{E}}[f(W)]=\widetilde{\mathcal{E}}[f(W)]\},$$
   \begin{equation}%\label{3}
   \nu\left(\lim_{n\rightarrow \infty}\frac{1}{\log n}\sum_{k\leq n}\frac{1}{k}f(W_{k}) = \widetilde{\mathbb{E}}\left[f(W)\right] \right)=1.
   \end{equation}
   \end{thm}
   %and $f$ in equation (\ref{3}) to be $f=f_{1}\circ f_{2}$, where $f_{2}:W(t) \mapsto W(1)$, $f_{1}\in \tilde{\mathcal{H}}$,
\fi

\begin{lem}\label{CLT} (Theorem 5.1 in Song \cite{Song2018}).
Let $\{X_i\}_{i=1}^{\infty}$ satisfy Assumption \ref{assump0} and \ref{assump}, then there exist a constant $\alpha\in(0,1)$ depending on $\beta$, and a constant $C_{\alpha,\beta}>0$ depending on $\alpha$, $\beta$ such that for any $n\ge 1$,
$$
\underset{|f|_{Lip}\leq 1}\sup \bigg| \mathbb{E}\left[ f\left( W_n \right) \right] -\widetilde{\mathbb{E}}\left[f(\xi)\right]\bigg| \leq C_{\alpha,\beta} \underset{1\leq i\leq n}\sup \bigg\{ \frac{\mathbb{E}[|X_i|^{2+\alpha}]}{\sigma_i^{2+\alpha}} \left( \frac{\sigma_i}{B_n} \right)^\alpha \bigg\},
$$
where $\xi$ is $G$-normal distribution under $\widetilde{\mathbb{E}}$ with the fixed $\beta$ and $\sqrt{\widetilde{\mathbb{E}}[\xi^2] }= \frac{2\beta}{1+\beta}$, $\sqrt{-\widetilde{\mathbb{E}}[-\xi^2] }= \frac{2 }{1+\beta}$.
\end{lem}

   \begin{thm}\label{thm1}
	Under Assumption \ref{assump0} and \ref{assump}, for the $\alpha$ in Lemma \ref{CLT}, if $\sup\limits_{i\geq 1}\mathbb{E} [|X_i|^{2+\alpha}] < \infty$,
	then for any $f \in \mathcal{H}$, we have
	\begin{equation}\label{3}
	\nu\left(\lim_{n\rightarrow \infty}\frac{1}{\log n}\sum_{k\leq n}\frac{1}{k}f\left(W_{k}\right) = \widetilde{\mathbb{E}}\left[f(\xi)\right] \right)=1,
	\end{equation}
where $\xi$ is $G$-normal distribution under $\widetilde{\mathbb{E}}$ with the fixed $\beta$ and $\sqrt{\widetilde{\mathbb{E}}[\xi^2] }= \frac{2\beta}{1+\beta}$, $\sqrt{-\widetilde{\mathbb{E}}[-\xi^2] }= \frac{2 }{1+\beta}$.
\end{thm}
We leave the proof of Theorem \ref{thm1} at the end of this section.
The following theorem is the main result of this paper.
   \begin{thm}[Almost Sure Central Limit Theorem]\label{thm3}
	Under Assumption \ref{assump0} and \ref{assump}, for the $\alpha$ in Lemma \ref{CLT}, if $\sup\limits_{i\geq 1}\mathbb{E} [|X_i|^{2+\alpha}] < \infty$,
	 we have
	\begin{equation}\label{eq2}
	\nu\left(\lim_{n\rightarrow \infty}\frac{1}{\log n}\sum_{k\leq n}\frac{1}{k}f\left(W_{k}\right) = \widetilde{\mathbb{E}}\left[f(\xi)\right],  \hbox{ for any } f \in \mathcal{H} \right)=1,
	\end{equation}
where $\xi$ is $G$-normal distribution under $\widetilde{\mathbb{E}}$ with the fixed $\beta$ and $\sqrt{\widetilde{\mathbb{E}}[\xi^2] }= \frac{2\beta}{1+\beta}$, $\sqrt{-\widetilde{\mathbb{E}}[-\xi^2] }= \frac{2 }{1+\beta}$.
\end{thm}
\Proof. Let $\{f_m\}_{m=1}^\infty$ be a countable dense subset of $\mathcal{H}$ under  $\|\cdot\|_\infty$.
For any $f\in \mathcal{H}$, there exists a subsequence of $\{f_m\}_{m=1}^\infty$, denoted by $\{f_{m_i}\}_{i=1}^\infty$, such that $f_{m_i}\rightarrow f$ under $\|\cdot\|_{\infty}$ as $i\to \infty$. Therefore, for any $\epsilon>0$, there exists an integer $i_1(\epsilon)\ge 1$ such that for all $i>i_1(\epsilon)$, we have
	$$
	\sup_{x\in\mathbb{R}} | f_{m_i}(x)-f(x) | < \epsilon.
	$$
%Thus,   we have
%	$$
%	| \widetilde{\mathbb{E}}[f_{m_i}(\xi)] - \widetilde{\mathbb{E}}[f(\xi)] |\le \widetilde{\mathbb{E}}[|f_{m_i}(\xi) - f(\xi)|]\le \sup_{x\in\mathbb{R}} | f_{m_i}(x)-f(x) | < \epsilon.
%	$$
Define
	$$
	A_m := \left\{ \omega: \lim_{n\rightarrow\infty} \frac{1}{\log n} \sum_{k\leq n}\frac1k f_m(W_k(\omega)) = \widetilde{\mathbb{E}}[f_m(\xi)] \right\}^c \hbox{ and } A:=\bigcup_{m=1}^\infty A_m.$$
According to Theorem 4.1, we have $\mathbb{V}(A_m) =0$, then $\mathbb{V}(A)\leq \sum_{m\geq 1 }\mathbb{V}(A_m)=0$, that is
\begin{eqnarray*}
\nu(A^c)=
\nu\left(\bigcap_{m=1}^{\infty}\left\{\omega:\ \lim_{n\rightarrow \infty}\frac{1}{\log n}\sum_{k\leq n}\frac{1}{k}f_m\left(W_{k}(\omega)\right) = \widetilde{\mathbb{E}}\left[f_m(\xi)\right]\right\} \right)=1.
\end{eqnarray*}
Fix any $\omega\in A^c$,  for all $m\ge1$, we have
$$\lim_{n\rightarrow \infty}\frac{1}{\log n}\sum_{k\leq n}\frac{1}{k}f_m\left(W_{k}(\omega)\right) = \widetilde{\mathbb{E}}\left[f_m(\xi)\right].$$
Hence, for $i>i_1(\epsilon)$, by the sub-additivity of $\widetilde{\mathbb{E}}$, we have
\begin{eqnarray*}
	&&\limsup_{n\to\infty}\left| \frac{1}{\log n} \sum_{k\leq n}\frac1k f(W_k(\omega)) - \widetilde{\mathbb{E}}[f(\xi)] \right|\\
 &\leq& \limsup_{n\to\infty}\left| \frac{1}{\log n} \sum_{k\leq n}\frac1k f(W_k(\omega)) - \frac{1}{\log n} \sum_{k\leq n}\frac1k f_{m_i}(W_k(\omega)) \right| \\
	& &+ \limsup_{n\to\infty}\left| \frac{1}{\log n} \sum_{k\leq n}\frac1k f_{m_i}(W_k(\omega)) - \widetilde{\mathbb{E}}[f_{m_i}(\xi)] \right| + \left| \widetilde{\mathbb{E}}[f_{m_i}(\xi)] - \widetilde{\mathbb{E}}[f(\xi)] \right| \\
 &\leq& \limsup_{n\to\infty} \frac{1}{\log n} \sum_{k\leq n}\frac1k \left| f(W_k(\omega)) - f_{m_i}(W_k(\omega)) \right|  +  \widetilde{\mathbb{E}}\left[\left|f_{m_i}(\xi)] - f(\xi) \right|\right]\\
	& \leq&  2\epsilon .
\end{eqnarray*}
Due to the arbitrariness of $\epsilon$,  we have
$$\lim_{n\to\infty} \frac{1}{\log n} \sum_{k\leq n}\frac1k f(W_k(\omega))= \widetilde{\mathbb{E}}[f(\xi)] .$$
Notice that $f$ is an arbitrary function in $\mathcal{H}$, so
$$A^c\subseteq \left\{\omega:\ \lim_{n\rightarrow \infty}\frac{1}{\log n}\sum_{k\leq n}\frac{1}{k}f\left(W_{k}(\omega)\right) = \widetilde{\mathbb{E}}\left[f(\xi)\right],  \hbox{ for any } f \in \mathcal{H} \right\}.$$
Therefore, equality (\ref{eq2}) is directly deduced from the monotonicity of $\nu$ and $\nu(A^c)=1$. The proof of Theorem \ref{thm3} is completed.
$\Box$

From Theorem \ref{thm3}, we can easily get the following corollary which is a version  of almost sure central limit theorem in Peng's sense.
\begin{coro}
Under Assumption \ref{assump0}, let $\{X_i\}_{i=1}^{\infty}$ be a sequence of independent and identically distributed random variables under a sub-linear expectation $\mathbb{E}$ with $\mathbb{E}[X_1]=\mathbb{E}[-X_1]=0$, $\bar{\sigma}=\sqrt{\mathbb{E}[X_1^2]}\ge \sqrt{-\mathbb{E}[-X_1^2]}=\underline{\sigma}>0$ and $\mathbb{E}[|X_1|^3]<\infty$. For $k\ge 1$, setting $S_k= \sum\limits_{i=1}^k X_i$, then
	\begin{equation*}
	\nu\left(\lim_{n\rightarrow \infty}\frac{1}{\log n}\sum_{k\leq n}\frac{1}{k}f\left(\frac{S_k}{\sqrt{k}}\right) = \widetilde{\mathbb{E}}\left[f(\xi)\right], \hbox{ for any } f \in \mathcal{H} \right)=1,
	\end{equation*}
where $\xi$ is $G$-normal distribution under $\widetilde{\mathbb{E}}$ with $G(a)=\frac{1}{2}(\overline{\sigma}^2a^{+}-\underline{\sigma}^2a^{-})$.
\end{coro}
  \begin{remark}
  If the expectation $\mathbb{E}$ in Theorem \ref{thm3} is the classical linear expectation, that is $\mathcal{P}=\{P\}$ is a singleton set,
   %, and    $\mathbb{E}[X]:=E_{P}[X]=\int X dP$,
   then we have $\mathcal{H}=C_{b,Lip}(\mathbb{R})$ and $\nu=P$. It means that Theorem \ref{thm3} is a fairly neat extension of classical almost sure central limit theorem in \cite{Brosamler}.
   \end{remark}
We will prove Theorem \ref{thm1} by using the strong law of large numbers (Theorem \ref{thm2}).
So we firstly give the following four lemmas.
For  all $ k,l,n \geq 1$, let $\xi_{k}:= f(W_k)- \mathbb{E}f(W_k)$, $Z_{l}:=\sum_{4^{l-1}\leq i < 4^{l}} \frac{\xi_{i}}{i}$ and $T_{n}:=\sum_{l=1}^{n} Z_{l}.$

\begin{lem}\label{lem}
   For all $f \in \mathcal{H}$, there exist two positive constants $M_{3}, M_{4}$ such that
   $$\mathbb{E}[\xi_{j}\xi_{k}] \leq M_{3}\cdot \left( \frac{j}{k}\right)^{\frac{1}{2}} + M_{4}\cdot\left(\frac{1}{k}\right)^{\frac{\alpha}{2}},$$
   where $j\leq k$ and $\alpha$ is the constant in  Lemma \ref{CLT}.
   \end{lem}
   \Proof. For $j\leq k$, let $r_j^k:= \frac{S_{k}-S_j}{B_k}$.
   By H\"{o}lder's  inequality (Lemma \ref{Holder}) and Rosenthal's inequality (Lemma \ref{Rosenthal}), we have
   \begin{align*}
   B_k \cdot \mathbb{E}\left[\bigg|\frac{S_{k}}{B_k} - r_{j}^{k}\bigg|\right] & %= \sqrt{k}\cdot \mathbb{E}\left[ \frac{|S_{j}|}{\sqrt{k}}\right]
   =\mathbb{E} \left[ |S_{j}|\right] \leq \sqrt{\mathbb{E}\left[ |S_{j}|^2\right]}
   %\leq \sqrt{C_{2}\cdot\overline{\sigma}^{2}\cdot j}
   \leq \sqrt{C_2 j \cdot \sup_{i\geq 1}\overline{\sigma}_i^2},
   \end{align*}
   where $C_{2}$ is the constant in Lemma \ref{Rosenthal}.
Then for all $ f \in C_{b,Lip}(\mathbb{R})$, $j\leq k$, by the definition of $B_k$, we have
   \begin{equation}\label{xijie}
   \mathbb{E}\bigg|f\left(W_k \right) -f(r_{j}^{k})\bigg| \leq \|f\|_{L} \mathbb{E}\left[\bigg|\frac{S_{k}}{B_k} - r_{j}^{k}\bigg|\right] \leq \|f\|_{L} \left( \frac{j}{k} \right)^{\frac{1}{2}}\sqrt{\frac{C_2 \sup\limits_{i\geq 1} \overline{\sigma}_i^2}{ \inf\limits_{i\geq 1} \underline{\sigma}_i^2}} ,
   \end{equation}
   where $\|f\|_{L}$ is the Lipschitz constant of $f$.
   %Because $\mathbb{E}[f(W_k)]\leq \mathbb{E}|f(W_k)-f(r_j^k)|+\mathbb{E}[f(r_j^k)]$, so $$\mathbb{E}[f(r_j^k)]-\mathbb{E}[f(W_k)]\geq -\mathbb{E}|f(r_j^k)-f(W_k)|,$$ in addition, $$\mathbb{E}[f(r_j^k)]-\mathbb{E}[f(W_k)]\leq \mathbb{E}|f(r_j^k)-f(W_k)|,$$ from above two, we have
   By the sub-additivity of $\mathbb{E}$, we have
   \begin{equation}\label{xiaren}
   \bigg|\mathbb{E}[f(r_j^k)] -\mathbb{E}\left[f\left(W_k \right)\right]\bigg|\leq \mathbb{E}\bigg|f(r_j^k)-f\left(W_k \right)\bigg|.
   \end{equation}

   Since $f$ is bounded, there exists a positive constant $M_{f}$ such that $|f|\leq M_{f}$. Due to  (\ref{xijie}) and (\ref{xiaren}), for all $j\leq k$, we have
   \begin{eqnarray*}
   &&\mathbb{E}[\xi_{j}\xi_{k}]\\
   &=&\mathbb{E}\left\{ \left[ f\left(W_j\right) -\mathbb{E}f \left(W_j\right)\right] \left[f\left(W_k\right)- \mathbb{E}f\left(W_k\right)\right]\right\} \\
   &=&\mathbb{E}\left\{\left[f\left(W_j\right) -\mathbb{E}f \left(W_j\right)\right] \left[f\left(W_k\right) -\mathbb{E}f\left(W_k\right)
   - [f(r_{j}^{k})-\mathbb{E}f(r_{j}^{k})] + [f(r_{j}^{k})-\mathbb{E}f(r_{j}^{k})]\right]\right\} \\
   &\leq &\mathbb{E}\left\{\left[f\left(W_j\right) -\mathbb{E}f\left(W_j \right)\right] \left[f\left(W_k \right)-f(r_{j}^{k})\right]\right\} \\ &\ \ \ & +\mathbb{E}\bigg\{\left[f\left(W_j \right)-\mathbb{E}f\left(W_j\right)\right] \left[\mathbb{E}f(r_{j}^{k}) -\mathbb{E}f\left(W_k \right)\right]\bigg\} \\
   &\ \ \ & +\mathbb{E}\bigg\{\left[f\left(W_j \right) -\mathbb{E}f\left(W_j \right)\right][f(r_{j}^{k})-\mathbb{E}f(r_{j}^{k})]\bigg\} \\
   &\leq &4 M_{f} \|f\|_{L} \left( \frac{j}{k}\right)^{\frac{1}{2}}  \sqrt{\frac{C_2 \sup\limits_{i\geq 1} \overline{\sigma}_i^2}{ \inf\limits_{i\geq 1} \underline{\sigma}_i^2}}  + \mathbb{E}\bigg\{\left[f\left(W_j \right) -\mathbb{E}f\left(W_j \right)\right][f(r_{j}^{k})-\mathbb{E}f(r_{j}^{k})]\bigg\} \\
   &\leq& 4M_{f}  \|f\|_{L} \left( \frac{j}{k}\right)^{\frac{1}{2}}\sqrt{\frac{C_2 \sup\limits_{i\geq 1} \overline{\sigma}_i^2}{ \inf\limits_{i\geq 1} \underline{\sigma}_i^2}} + \mathbb{E}\bigg\{\left[f\left(W_j \right) -\mathbb{E}f\left(W_j \right)\right]^{+}\mathbb{E}[f(r_{j}^{k})-\mathbb{E}f(r_{j}^{k})] \\ &\ \ \ &+\left[f\left(W_j \right) -\mathbb{E}f\left(W_j \right)\right]^{-}\mathbb{E}[-f(r_{j}^{k})+\mathbb{E}f(r_{j}^{k})]\bigg\}\\
     &=& 4M_{f} \|f\|_{L}  \left( \frac{j}{k}\right)^{\frac{1}{2}}\sqrt{\frac{C_2 \sup\limits_{i\geq 1} \overline{\sigma}_i^2}{ \inf\limits_{i\geq 1} \underline{\sigma}_i^2}}  + \mathbb{E}\bigg\{\left[f\left(W_j \right) -\mathbb{E}f\left(W_j \right)\right]^{-}\mathbb{E}[-f(r_{j}^{k})+\mathbb{E}f(r_{j}^{k})]\bigg\}\\
       &\leq& 4M_{f} \|f\|_{L}  \left( \frac{j}{k}\right)^{\frac{1}{2}}\sqrt{\frac{C_2 \sup\limits_{i\geq 1} \overline{\sigma}_i^2}{ \inf\limits_{i\geq 1} \underline{\sigma}_i^2}}  + 2M_{f} [\mathbb{E}[f(r_{j}^{k})]+\mathbb{E}[-f(r_{j}^{k})]]
   \end{eqnarray*}
  where the penultimate inequality is deduced from $f(r_j^k)$ being independent of $f\left(W_j\right)$ under $\mathbb{E}$. It follows from (\ref{xiaren}) that
   \begin{align*}
   \mathbb{E}[\xi_{j}\xi_{k}] &\leq 4M_{f}\|f\|_{L}  \left( \frac{j}{k}\right)^{\frac{1}{2}}  \sqrt{\frac{C_2 \sup\limits_{i\geq 1} \overline{\sigma}_i^2}{ \inf\limits_{i\geq 1} \underline{\sigma}_i^2}} + 2M_{f} \bigg\{\mathbb{E}[f(r_{j}^{k})] -\mathbb{E}\left[f\left(W_k \right)\right] \\ &\ \ \ +\mathbb{E}\left[f\left(W_k \right)\right] +\mathbb{E}\left[-f\left(W_k \right)\right] -\mathbb{E}\left[-f\left(W_k \right)\right] +\mathbb{E}[-f(r_{j}^{k})]\bigg\} \\
   &\leq 4M_{f}\|f\|_{L}  \left( \frac{j}{k}\right)^{\frac{1}{2}}\sqrt{\frac{C_2 \sup\limits_{i\geq 1} \overline{\sigma}_i^2}{ \inf\limits_{i\geq 1} \underline{\sigma}_i^2}}  + 2M_{f} \bigg\{\mathbb{E}\left[|f(r_{j}^{k})-f\left(W_k \right)|\right] \\ &\ \ \ +\mathbb{E}\left[f\left(W_k \right)\right]+\mathbb{E}\left[-f\left(W_k \right)\right] +\mathbb{E}\left[|f\left(W_k \right)-f(r_{j}^{k})|\right]\bigg\}.
   \end{align*}
Moreover, by (\ref{xijie}) and $\widetilde{\mathbb{E}}\left[f(\xi)\right]=-\widetilde{\mathbb{E}}\left[-f(\xi)\right]$, we have
   \begin{eqnarray*}
   &&\mathbb{E}[\xi_{j}\xi_{k}]\\
   &\leq& 4M_{f} \|f\|_{L}  \left( \frac{j}{k}\right)^{\frac{1}{2}}\sqrt{\frac{C_2 \sup\limits_{i\geq 1} \overline{\sigma}_i^2}{ \inf\limits_{i\geq 1} \underline{\sigma}_i^2}}  \\ && + 2M_{f} \left\{  2 \|f\|_{L} \left( \frac{j}{k}\right)^{\frac{1}{2}} \sqrt{\frac{C_2 \sup\limits_{i\geq 1} \overline{\sigma}_i^2}{ \inf\limits_{i\geq 1} \underline{\sigma}_i^2}} +\mathbb{E}\left[f\left(W_k \right)\right] +\mathbb{E}\left[-f\left(W_k \right)\right] \right\} \\
   %&=8M_{f} \overline{\sigma}\|f\|_{L}  \left( \frac{j}{k}\right)^{\frac{1}{2}} + 2M_{f} \{\mathbb{E}[f(W_{k})] -\mathcal{E}[f(W_{k})] \} \\
   &\leq& 8M_{f}\|f\|_{L}  \left( \frac{j}{k}\right)^{\frac{1}{2}} \sqrt{\frac{C_2 \sup\limits_{i\geq 1} \overline{\sigma}_i^2}{ \inf\limits_{i\geq 1} \underline{\sigma}_i^2}} \\
    && + 2M_{f} \left\{ \mathbb{E}\left[f\left(W_k \right)\right]-\widetilde{\mathbb{E}}\left[f(\xi)\right]+\widetilde{\mathbb{E}}\left[f(\xi)\right]-\widetilde{\mathbb{E}}\left[-f(\xi)\right] +\widetilde{\mathbb{E}}\left[-f(\xi)\right]                        +\mathbb{E}\left[-f\left(W_k \right)\right] \right\} \\
&\leq&8M_{f}\|f\|_{L}  \left( \frac{j}{k}\right)^{\frac{1}{2}} \sqrt{\frac{C_2 \sup\limits_{i\geq 1} \overline{\sigma}_i^2}{ \inf\limits_{i\geq 1} \underline{\sigma}_i^2}}  + 4M_{f} C_{\alpha,\beta}\sup_{1\le i\le k}\frac{\mathbb{E}[|X_i|^{2+\alpha}]}{\sigma_i^{2+\alpha} \left(\frac{\sqrt{\sum_{i=1}^{k}\sigma_i^2}}{\sigma_i}\right)^{\alpha}} \\
   &\leq& 8M_{f}\|f\|_{L}  \left( \frac{j}{k}\right)^{\frac{1}{2}} \sqrt{\frac{C_2 \sup\limits_{i\geq 1} \overline{\sigma}_i^2}{ \inf\limits_{i\geq 1} \underline{\sigma}_i^2}}  + 4M_{f} C_{\alpha,\beta} \left(\frac{1}{k}\right)^{\frac{\alpha}{2}} \frac{\sup_{i\geq 1}\mathbb{E}[|X_i|^{2+\alpha}]}{(\inf_{i\geq 1}\sigma_i^2)^{\frac\alpha2+1}} \\
   &:=&M_{3}\cdot \left( \frac{j}{k}\right)^{\frac{1}{2}} + M_{4}\cdot\left(\frac{1}{k}\right)^{\frac{\alpha}{2}} ,
   \end{eqnarray*}
   where the penultimate inequality is deduced from Lemma \ref{CLT}. %Actually, for all $ f \in \mathcal{H}$
   %\begin{align*}
   %\mathbb{E}\left[f\left(W_k \right)\right] - \mathcal{E}\left[f\left(W_k \right)\right] &= \mathbb{E}\left[f\left(W_k \right)\right] -\widetilde{\mathbb{E}}[f(\xi)] +\widetilde{\mathcal{E}}[f(\xi)] -\mathcal{E}\left[f\left(\frac{S_{k}}{\sqrt{k}}\right)\right] \\
   %& \leq 2 C\left(\frac{1}{k}\right)^{\frac{\alpha}{2}}(1+\mathbb{E}[|X_{1}|^{\alpha}+|X_{1}|^{2+\alpha}]),
   %\end{align*}
   %where $C$ is the constant in Theorem A.1.
   Therefore, the proof of this lemma is completed. $\Box$

   \begin{lem}\label{lem'}
  For all $ f \in \mathcal{H}$, there exist two positive constants $M_{3}, M_{4}$ such that
   $$\mathbb{E}Z_{l}Z_{m} \leq 9M_{3}\cdot 2^{l-m} + 9M_{4}\cdot 2^{\alpha(1-m)},$$
   where $l<m$, $\alpha$, $M_3$ and $M_4$ are the constants in Lemma \ref{lem}.
   \end{lem}
   \Proof. For $l<m$, it follows from Lemma \ref{lem} that
   \begin{align*}
   \mathbb{E}Z_{l}Z_{m} &=\mathbb{E}\left[\sum_{4^{l-1}\leq i < 4^{l}} \frac{\xi_{i}}{i} \sum_{4^{m-1}\leq j < 4^{m}} \frac{\xi_{j}}{j} \right] \\
   &\leq \sum_{i,j} \frac{\mathbb{E}[\xi_{i}\xi_{j}]}{ij} \\
   &\leq M_{3}\sum_{i,j}\frac{\left( \frac{i}{j}\right)^{\frac{1}{2}}}{ij} + M_{4}\sum_{i,j}\frac{\left( \frac{1}{j}\right)^{\frac{\alpha}{2}}}{ij} \\
   &= M_{3}\sum_{i,j}\frac{1}{i^{\frac{1}{2}}j^{\frac{3}{2}}} + M_{4}\sum_{i,j}\frac{1}{ij^{1+\frac{\alpha}{2}}} \\
   &\leq 9M_{3} \frac{4^{l-1}4^{m-1}}{2^{l-1}2^{3(m-1)}} + 9M_{4} \frac{4^{m-1}}{4^{(m-1)(1+\frac{\alpha}{2})}} \\
   &= 9M_{3} 2^{l-m} + 9M_{4} 2^{\alpha(1-m)}.
   \end{align*}
   Thus, we complete the proof of Lemma \ref{lem'}. $\Box$

   \begin{lem}\label{lem6}
   For any $n_1\le n_2$, there exists a positive constant $M_{5}$ such that
   $$\sum_{l=n_{1}}^{n_{2}}\mathbb{E}[Z_{l}^2]\leq M_{5}\cdot(n_{2}-n_{1}+1).$$
   \end{lem}
   \Proof. By Lemma \ref{lem} we have
   \begin{align*}
   \sum_{l=n_{1}}^{n_{2}}\mathbb{E}[Z_{l}^2] &= \sum_{l=n_{1}}^{n_{2}} \mathbb{E} \left[ \sum_{4^{l-1}\leq i < 4^{l}} \frac{\xi_{i}}{i} \right]^2 = \sum_{l=n_{1}}^{n_{2}} \mathbb{E} \left[ \sum_{4^{l-1}\leq i < 4^{l}} \frac{\xi_{i}^2}{i^2} + 2\sum_{4^{l-1}\leq i <j <4^{l}} \frac{\xi_{i}\xi_{j}}{ij}\right] \\
   &\leq \sum_{j=4^{n_{1}-1}}^{4^{n_{2}}-1}\mathbb{E}\frac{\xi_{j}^2}{j^2}+ 2\sum_{l=n_{1}}^{n_{2}}\sum_{4^{l-1}\leq i <j <4^{l}}\frac{M_{3}( \frac{i}{j})^{\frac{1}{2}} + M_{4}\cdot\left(\frac{1}{j}\right)^{\frac{\alpha}{2}}}{ij} \\
   &\leq \sum_{j=4^{n_{1}-1}}^{4^{n_{2}}-1}\frac{4M_{f}^2}{j^2} + 2M_{3}\sum_{l=n_{1}}^{n_{2}}\sum_{4^{l-1}\leq i <j <4^{l}}\frac{1}{i^{\frac{1}{2}}j^{\frac{3}{2}}}  + 2M_{4}\sum_{l=n_{1}}^{n_{2}}\sum_{4^{l-1}\leq i <j <4^{l}}\frac{1}{ij^{1+\frac{\alpha}{2}}}\\
   &\leq \sum_{j=1}^{\infty}\frac{4M_{f}^2}{j^2} + 2M_{3}\sum_{l=n_{1}}^{n_{2}}9 + 2M_{4}\sum_{l=n_{1}}^{n_{2}}9\cdot 2^{\alpha(1-l)} \\
   &\leq \sum_{j=1}^{\infty}\frac{4M_{f}^2}{j^2} + 4(M_{3}+M_{4})\sum_{l=n_{1}}^{n_{2}}9\\
   &\leq M_{5}\cdot(n_{2}-n_{1}+1).
   \end{align*}
   Hence, the proof of Lemma \ref{lem6} is completed. $\Box$

   \begin{lem}\label{lem7}
For all $n\ge 1$, there exists a positive constant $M_{6}$ such that
   $$\mathbb{E}[T_{n}^2]\leq M_{6}\cdot n.$$
   \end{lem}
   \Proof. According to Lemma \ref{lem'} and Lemma \ref{lem6}, we have
   \begin{align*}
   \mathbb{E}[T_{n}^2] &=\mathbb{E}\left[ \sum_{l=1}^{n} Z_{l}^2 + 2 \sum_{1\leq l<m\leq n} Z_{l}Z_{m}\right]
   \leq \sum_{l=1}^{n}\mathbb{E}Z_{l}^2 +2 \sum_{1\leq l<m\leq n} \mathbb{E}[Z_{l}Z_{m}] \\
   &\leq M_{5}n + 18M_{3} \sum_{1\leq l<m\leq n} 2^{l-m} + 18M_{4} \sum_{1\leq l<m\leq n} 2^{\alpha(1-m)} \\
   &\leq M_{5}n + 18M_{3}(n-1) \sum_{m=1}^{n-1} \frac{1}{2^{m}} + 18M_{4}(n-1)2^{\alpha}\sum_{m=1}^{n}\frac{1}{2^{\alpha m}} \\
   &\leq  M_{5}n + 18M_{3}(n-1)\sum_{m=1}^{\infty} \frac{1}{2^{m}} + 18M_{4}(n-1)2^{\alpha}\sum_{m=1}^{\infty}\frac{1}{2^{\alpha m}}
   =M_{6}n.
   \end{align*}
   This completes the proof of Lemma \ref{lem7}. $\Box$

   Now, we turn to the proof of Theorem \ref{thm1}.

   \Proof. [Proof of Theorem \ref{thm1}]
   For each $f \in \mathcal{H}$, if the following equality holds
   \begin{equation}\label{1}
   \nu\left( \limsup_{n\rightarrow \infty}\frac{1}{\log n}\sum_{k\leq n}\frac{1}{k}f\left( W_k \right) \leq \widetilde{\mathbb{E}}\left[f(\xi)\right] \right)=1,
   \end{equation}
   then
   \begin{equation*}
   \nu\left( \liminf_{n\rightarrow \infty}\frac{1}{\log n}\sum_{k\leq n}\frac{1}{k}f\left( W_k \right) \geq -\widetilde{\mathbb{E}}\left[-f(\xi)\right] \right)=1,
   \end{equation*}
    since $-f\in\mathcal{H}$.
   Therefore (\ref{3}) holds, that is
   \begin{equation*}
   \nu\left( \lim_{n\rightarrow \infty}\frac{1}{\log n}\sum_{k\leq n}\frac{1}{k}f\left( W_k \right) = \widetilde{\mathbb{E}}[f(\xi)] \right)=1.
   \end{equation*}
   So we only need to prove (\ref{1}).

   By Lemma \ref{CLT}, we have
   $$\lim_{k\rightarrow\infty}\mathbb{E}\left[f\left( W_k \right)\right]= \widetilde{\mathbb{E}}[f(\xi)].$$
 Hence, for any $ \epsilon >0$, there exists $ K>0$, such that for all $ k>K$,
   $$\widetilde{\mathbb{E}}[f(\xi)]-\epsilon \leq \mathbb{E}\left[f\left( W_k \right)\right] \leq \widetilde{\mathbb{E}}f(\xi)+\epsilon,$$
   and
   $$\frac{1}{\log n}\sum_{K\leq k\leq n}\frac{1}{k}(\widetilde{\mathbb{E}}[f(\xi)]-\epsilon) \leq \frac{1}{\log n}\sum_{K\leq k\leq n}\frac{1}{k}\mathbb{E}\left[f\left( W_k \right)\right] \leq \frac{1}{\log n}\sum_{K\leq k\leq n}\frac{1}{k}(\widetilde{\mathbb{E}}[f(\xi)]+\epsilon).$$
   Note that
   $$\lim_{n\rightarrow\infty}\frac{1}{\log n}\sum_{K\leq k\leq n}\frac{1}{k}=1,$$
   we have
   $$\lim_{n\rightarrow\infty}\frac{1}{\log n}\sum_{K\leq k\leq n}\frac{1}{k}\mathbb{E}\left[f\left(W_k \right)\right]= \widetilde{\mathbb{E}}f(\xi),$$
   and
   $$\lim_{n\rightarrow\infty}\frac{1}{\log n}\sum_{k\leq n}\frac{1}{k}\mathbb{E}f\left( W_k \right) = \lim_{n\rightarrow\infty}\frac{1}{\log n}\sum_{k\leq K}\frac{1}{k}\mathbb{E}f\left(W_k \right) +\lim_{n\rightarrow\infty}\frac{1}{\log n}\sum_{K\leq k\leq n}\frac{1}{k}\mathbb{E}f\left(W_k \right)=\widetilde{\mathbb{E}}f(\xi).$$
   Therefore, (\ref{1}) is equivalent to
   $$\nu\left(\limsup_{n\rightarrow\infty}\frac{1}{\log n}\sum_{k\leq n}\frac{1}{k}f\left( W_k \right)\leq  \lim_{n\rightarrow\infty}\frac{1}{\log n}\sum_{k\leq n}\frac{1}{k}\mathbb{E}f\left( W_k \right) \right)=1,$$
   which is equivalent to
   $$\nu\left(\limsup_{n\rightarrow\infty}\frac{1}{\log n}\sum_{k\leq n}\frac{1}{k}\left(f\left( W_k \right) -\mathbb{E}\left[f\left( W_k \right)\right]\right)\leq 0 \right)=1.$$
   Let $\xi_{k}:=f\left( W_k \right)-\mathbb{E}f\left( W_k \right)$, then (\ref{1}) is equivalent to
   $$\nu\left(\limsup_{n\rightarrow\infty}\frac{1}{\log n}\sum_{k\leq n}\frac{1}{k}\xi_{k}\leq 0 \right)=1,$$
   which is equivalent to
   $$\nu\left(\limsup_{n\rightarrow\infty}\frac{1}{\log_{4}n}\sum_{k\leq n}\frac{1}{k}\xi_{k}\leq 0 \right)=1.$$
   Let $N:=\log_{4}n$, then (\ref{1}) is equivalent to
   $$\nu\left( \limsup_{N\rightarrow \infty} \frac{1}{N}\sum_{l=1}^{N} Z_{l}\leq 0 \right)=1, $$
   where $Z_{l}:=\sum_{4^{l-1}\leq k < 4^{l}} \frac{\xi_{k}}{k}$, for all $l \geq 1.$
  For all $n\ge 1$, let $T_{n}:=\sum_{l=1}^{n}Z_{l}$, then (\ref{1}) is equivalent to
   $$\nu\left( \limsup_{n\rightarrow \infty} \frac{T_{n}}{n}\leq 0 \right)=1. $$
   It is easy to check that
   $$\mathbb{E}[T_{n}]\leq \sum_{l=1}^{n} \mathbb{E}[Z_{l}] \leq \sum_{l=1}^n\sum_{4^{l-1}\le k<4^l} \frac{\mathbb{E}[\xi_{k}]}{k} = \sum_{l=1}^n\sum_{4^{l-1}\le k<4^l} \frac{\mathbb{E}[f(W_k)-\mathbb{E}f(W_k )]}{k} =0.$$
   By Lemma \ref{lem6} and Lemma \ref{lem7}, we have $\sum_{l=n_{1}}^{n_{2}}\mathbb{E}(Z_{l}^2)\leq M_{5}\cdot (n_{2}-n_{1}+1)$ and $\mathbb{E}[T_{n}^{2}]\leq M_6 \cdot n$,  for all $n_2\ge n_1$ and $n \geq 1$. Therefore, equality (\ref{1}) can be deduced from Theorem \ref{thm2}.
   The proof of Theorem \ref{thm1} is completed. $\Box$
%=================================================================================================================

\hspace{-0.6cm}\textbf{Acknowledgements}\\The work is supported by the National Natural Science Foundation of China (Grant Nos. 11601280 and 11701331 and 11471190) and the Natural Science Foundation of Shandong Province of China
(Grant Nos. ZR2017AQ007 and ZR2016AQ11 and ZR2016AQ13).

\end{document}